\documentclass[aap,MSNbibl,seceqn,dvips]{arximspdf}
\usepackage{graphicx}

\doi{10.1214/13-AAP1001} 
\volume{25}
\issue{1}
\pubyear{2015}
\firstpage{406}
\lastpage{428}

\makeatletter

\newcommand{\rrVert}{\Vert}
\newcommand{\llVert}{\Vert}
\def\cal{\mathcal}
\newcommand{\eqref}[1]{(\ref{#1})}
\newtheorem{theorem}{Theorem}[section]
\newtheorem{lemma}[theorem]{Lemma}
\newproclaim{assumption}[theorem]{Assumption}
\newproclaim{example}{Example}
\newproclaim{definition}[theorem]{Definition}
\newtheorem{corollary}[theorem]{Corollary}
\newproclaim{remark}[theorem]{Remark}

\newcommand{\la}{\lambda}
\newcommand{\eps}{\varepsilon}
\newcommand{\X}{\Xi}
\newcommand{\Ps}{\Psi}

\newcommand{\Om}{\Omega}
\newcommand{\Z}{\mathbb Z}
\newcommand{\R}{{\mathbb R}}
\newcommand{\N}{{\mathbb N}}
\newcommand{\PP}{{\mathbb P}}
\newcommand{\D}{{\mathbb D}}
\newcommand{\calE}{{\cal E}}
\newcommand{\calF}{{\cal F}}
\newcommand{\calG}{{\cal G}}
\newcommand{\calI}{{\cal I}}
\newcommand{\calJ}{{\cal J}}

\newcommand{\calP}{{\cal P}}
\newcommand{\calV}{{\cal V}}

\newcommand{\cd}{\cdot}
\makeatother

\begin{document}
\begin{frontmatter}

\title{Necessary condition for null controllability
in many-server heavy traffic}
\runtitle{Necessary condition for null controllability}

\begin{aug}
\author[A]{\fnms{Gennady} \snm{Shaikhet}\corref{}\thanksref{T1}\ead[label=e1]{gennady@math.carleton.ca}}
\runauthor{G. Shaikhet}
\thankstext{T1}{Supported in part by the NSERC Grant N: RGPIN 402292-2011.}
\affiliation{Carleton University}
\address[A]{Department of Mathematics and Statistics\\
Carleton University\\
Ottawa, Ontario K1S 5B6\\
Canada\\
\printead{e1}}
\end{aug}

\received{\smonth{9} \syear{2013}}
\revised{\smonth{12} \syear{2013}}

%
\begin{abstract}
Throughput sub-optimality (TSO), introduced in Atar and Shaik\-het
[\textit{Ann. Appl. Probab.} \textbf{19} (2009) 521--555] for static fluid
models of parallel queueing networks, corresponds to the existence of a
resource allocation, under which the total service rate becomes greater
than the total arrival rate. As shown in Atar, Mandelbaum and Shaikhet
[\textit{Ann. Appl. Probab.} \textbf{16} (2006) 1764--1804] and Atar and
Shaikhet (2009), in the many server Halfin--Whitt regime, TSO implies
null controllability (NC), the existence of a routing policy under
which, for every finite $T$, the measure of the set of times prior to
$T$, at which at least one customer is in the buffer, converges to zero
in probability at the scaling limit. The present paper investigates the
question whether the converse relation is also true and TSO is both
sufficient and necessary for the NC behavior.

In what follows we do get the affirmation for systems with either two
customer classes (and possibly more service pools) or vice-versa and
state a condition on the underlying static fluid model that allows the
extension of the result to general structures.
\end{abstract}

%
\begin{keyword}[class=AMS]
\kwd{60K25}
\kwd{68M20}
\kwd{90B22}
\kwd{90B36}
\kwd{60F05}
\end{keyword}
\begin{keyword}
\kwd{Multiclass queueing systems}
\kwd{heavy traffic}
\kwd{scheduling and routing}
\kwd{throughput optimality}
\kwd{null controllability}
\kwd{buffer-station graph}
\kwd{simple paths}
\end{keyword}

\end{frontmatter}

\section{Introduction}

In this paper we consider many-server parallel queueing
networks in heavy traffic regime. Despite the criticality, as shown in
\cite{AMS,AS}, there may exist a scheduling rule, with high probability
maintaining the system without waiting customers, for ``most of the
time.'' Called \textit{null controllability}, such unusual phenomena
occurs under the \textit{throughput sub-optimality} of the underlying
(critically loaded in a standard sense), fluid model. In the current
work we try to understand if the effect can still be achieved when the
underlying fluid is throughput optimal, and conclude that it is not
possible and throughput sub-optimality is indeed required.

\begin{figure}

\includegraphics{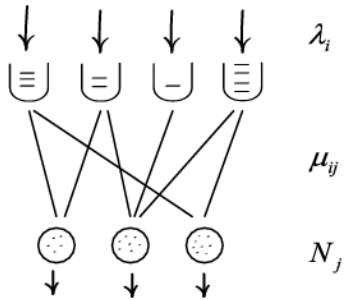}

\caption{A queueing model with four customer classes and three
service pools.}\label{fig1}
\end{figure}

Our model of interest consists of multiple customer classes, indexed by
$\calI$, and several service pools, indexed by $\calJ$, each consisting
of many i.i.d. exponential servers. The servers rates depend on both
the station and the class.
A system administrator dynamically controls all scheduling and routing
in the system; see Figure~\ref{fig1}. The model is considered in the heavy
traffic parametric regime, first proposed by Halfin and Whitt \cite
{halwhi}, in which
the number of servers at each station and the arrival rates
grow without bound, proportionally to some $n\uparrow\infty$.

Typically, when analyzing such systems, one looks at the underlying
\textit{static fluid model}, obtained in the law of large numbers limit
of the processes involved. According to \cite{harlop,wil}, the
so-called \textit{static fluid allocation problem} (see Section~\ref
{SFM} for the details) should then be formulated to determine whether
the model (hence, the original system) is under, over or critically
loaded; the latter being the proper foundation for the heavy traffic
analysis. What one gets is a deterministic matrix~$\xi^*$, where for
$(i,j)\in\calI\times\calJ$, the entry $\xi^*_{ij}$ represents the
fraction of station-$j$ work dedicated to class-$i$ customers on the
fluid level. Consequently, the original network is called critically
loaded if the optimized fluid takes $100 \%$ of system capacity; that
is, $\sum_i \xi_{ij}^*=1$ for each $j\in\calJ$. The class-station
pairs $(i,j)\in\calI\times\calJ$, along which the service is
possible, are called \textit{activities}. The activities $(i,j)$ for
which $\xi^*_{ij}>0$ (resp., $\xi_{ij}=0$)
are regarded as \textit{basic} (resp., \textit{nonbasic}). In both
\cite{harlop,wil} the set $\xi^*$ was assumed to be unique as well as
to satisfy the \textit{complete resource pooling} condition, requiring
all vertices in the class-station graph to be connected via basic
activities. Under the uniqueness assumption, the latter was shown to be
equivalent for the graph of basic activities being a tree. The above
set of conditions on the underlying fluid model has become standard for
considerable amount of works in the conventional (e.g., \cite
{harlop,wil,mansto}, etc.) and Halfin--Whitt (e.g., \cite
{atar1,atar2,atamanrei}, etc.) heavy traffic regimes, as well as in
the more
recent, nondegenerate-slowdown (NDS) regime \cite{atagur}.

With the static fluid model set, an attempt is then made
to prove that appropriately scaled (Halfin--Whitt regime) fluctuations
of the queueing model about
the fluid model converge to a diffusion. Assuming no use of nonbasic
activities, the pioneering papers \cite{atar1,atar2,atamanrei} were
able to represent the scaled system dynamics as a controlled diffusion
with a \textit{drift control}, thus being able to determine
asymptotically optimal scheduling policies for the fairly large class
of operational costs. It is a general understanding (see Theorem~2.1 in
\cite{AMS} and Theorem~\ref{12121212} in this paper) that by including the
nonbasic activities one gets additional controls, this time \textit
{singular controls}, but the augmented controlled diffusion still
remains to be fully analyzed.


Yet, some partial results had been obtained. One of them, called
\textit{null controllability} will be the focus of our attention. In
particular, in \cite{AMS}, Theorems 2.3, 2.4, it was shown that in
the presence of nonbasic activities, some models are prone to a fairly
unusual effect when a critically loaded system starts to behave like an
underloaded system. More exactly, one can construct a policy, under
which for any given $0<\eps<T<\infty$
all queues in the system are kept empty
on the time interval $[\eps,T]$, with probability approaching one (a
finiteness of the interval in \cite{AMS} is crucial, and was later
supported by the results from \cite{stotez} indicating that the
phenomenon is not possible in the long run). It is also worth noting
that null controllability seems to be the feature of the Halfin--Whitt
regime only---by its nature, it cannot happen in the conventional
single server asymptotics---and the \textit{conventional-like} NDS
regime is not expected to have it either.

The results of \cite{AMS} were generalized and better explained in
\cite
{AS}, attributing the null controllability to what was called \textit
{throughput sub-optimality} of the underlying static fluid model, a
situation, when (static) resources can be rearranged so that the total
service rate becomes greater than the total arrival rate; see
Section~\ref{TO} for the exact definition. Throughput sub-optimality, it
appears, may occur in wider class of fluid models and, even when the
null controllability is impossible, can result in its weaker (though
still efficient) version. Namely (Theorem~1 of \cite{AS}), for every
finite $T$, the measure of the set of times prior to $T$, at which at
least one customer is in the buffer, converges to zero in probability
at the scaling limit.

This brings us to the main objective, to understand the converse
relation between throughput sub-optimality and (weak) null
controllability (Theorem~\ref{main_thm}). We show that the desired
property is rooted in a simply formulated
deterministic result (Theorem~\ref{elimin}) stating that a throughput
optimal static fluid model does not become sub-optimal if its fluid
amounts are modified along the so-called zero paths, simple paths $p$
from \cite{AS} with signed weight $\mu(p)=0$. This gives a new
interesting perspective on zero paths, normally not associated with
abrupt changes of fluid material; in contrast with ``unwelcome''
positive paths $\mu(p)>0$ that increase the material, or
negative paths $\mu(p)<0$, the existence of which, as shown in Theorem~2 of \cite{AS}, is equivalent to throughput sub-optimality. Both
Theorem~\ref{elimin} and its dynamic version Lemma~\ref{estimate_main}
are proven for systems with either two customer classes (and possibly
more service pools) or vice-versa. Although the full version of Theorem~\ref{elimin} still remains unresolved, its simplistic nature [checkable
relations \eqref{AA4}--\eqref{AA5}] allows us to partially generalize
the results for arbitrary $I$ and $J$ (Theorem~\ref{main_thm_gen}).

%

The organization of the paper is rather straightforward, with the main
result (Theorem~\ref{main_thm}) followed by its proof (Section~\ref
{proof}). Sections~\ref{OM}--\ref{SFM} provide all the
prerequisites, while Section~\ref{road} is the roadmap for the
proof. After that, Theorem~\ref{main_thm_gen} of Section~\ref{gs}
discusses possible extensions of our findings.

\textit{Notation}. For a positive integer $d$ and $x\in\R^d$, let
$\|x\|=\sum_{i=1}^d|x_i|$. For $v,u\in\R^d$ let $v\cd
u=\sum_{i=1}^du_iv_i$. The symbols $e_i$ denote the unit coordinate
vectors and $e=(1,\ldots,1)$. The dimension of $e$ may change from
one expression to another. Thus, for $x=(x_1,\ldots, x_d) \in\R^d$, we
have $e\cdot x = \sum_{i=1}^d x_i$. For an event $A$ we use
$ 1\{A\}$ for the indicator of $A$. Denote by $\D(\R^d)$ the space of all
cadlag functions (i.e., right continuous and having finite left
limits) from $\R_+$ to $\R^d$. Denote
$\vert X\vert^*_t=\sup_{0\le u\le t}\vert X(u)\vert$ for $X\in\D
(\R)$,
$\Vert X\Vert^*_t=\sup_{0\le u\le t}\Vert X(u)\Vert$ for $X\in\D
(\R^d)$ and
$f(t:s) = f(t)-f(s)$.

%
%
%

\section{The model and the main result} \label{ch2}

\subsection{Original model}\label{OM}

The setting is standard; see, for example, \cite
{atar2,atamanrei,AMS,AS}. A~complete probability space $(\Om,\calF,\PP)$ is given, supporting
all stochastic processes defined below. There is a sequence of systems
indexed by $n\in\N$, each having $I$ customer classes and $J$ service
stations. Station $j$ has $N_j^n$ identical servers. The classes are
labeled as $1,\ldots,I$ and the stations as $I+1,\ldots,I+J$. We set
$ \calI=\{1,\ldots,I\}, \calJ=\{I+1,\ldots,I+J\}$. The arrival and
service processes, all mutually independent, are denoted by $\{A^n_i,
i\in\calI\}$ and $\{S_{ij}^n, (i,j)\in\calI\times\calJ\}$. Each
$A_i^n$ is a renewal process whose inter-arrival times have finite
second moment and an inverse mean (or rate) equal to $\la_i^n>0$. Each
service process $S_{ij}^n$ is a Poisson process with rate $\mu
^n_{ij}\ge0$. We also allow a possibility for
$\mu_{ij}^n=0$, in which case we say that class-$i$ customers cannot
be served at station $j$.

Denote the set of all
class-station pairs by $\calE:=\calI\times\calJ$, let
$
\calE_\mathrm{ a}=\{(i,j)\in\calI\times\calJ\dvtx\mu_{ij}^n>0\}$,
and, throughout, assume that $\calE_\mathrm{ a}$ does not depend on
$n$. A class-station pair $(i,j)\in\calE_\mathrm{ a}$ is said to be an
\textit{activity}. The set of class-station pairs that are not
activities is
denoted by $\calE_\mathrm{ a}^c\equiv\calE\setminus\calE_\mathrm{ a}$.

The number of service completions of class-$i$ customers by all
servers of station $j$ by time $t$ is therefore (see, e.g., \cite
{atar2,atamanrei,AMS,AS}), given by
$S_{ij}^n (\int_0^t\Ps_{ij}^n(s)\,ds )$, where for every
$(i,j)\in
\calE_\mathrm{ a}$, we
denote by $\Ps_{ij}^n(t)$ the number of class-$i$ customers being
served in station $j$ at time $t$. Denote by $X^n_i(t)$ the number of
class-$i$ customers in the
system at time $t$. By definition,
%
\begin{eqnarray}
\label{23} X_i^n(t)&=&X_i^{n}(0)+A_i^n(t)
-\sum_{j\in\calJ}S_{ij}^n \biggl(\int
_0^t\Ps_{ij}^n(s)\,ds
\biggr),\qquad i\in\calI;
\\
\label{20} \qquad\sum_{j\in\calJ}\Ps^n_{ij}(t)
&\le& X_i^n(t), \qquad i\in\calI;\qquad  \sum
_{i\in\calI}\Ps^n_{ij}(t)\le
N_j^n,\qquad j\in\calJ.
\end{eqnarray}
The processes $\Ps^n=(\Ps^n_{ij})_{(i,j)\in\calI\times\calJ}$ are
regarded as
\textit{scheduling control policy} (SCP) and assumed to be
right-continuous, taking values in $\Z_+$. Thus
%
\begin{equation}
\label{40} \Ps_{ij}^n(t)\ge0,\qquad (i,j)\in
\calE_\mathrm{ a};\qquad \Ps_{ij}^n(t)=0,\qquad (i,j) \in
\calE_\mathrm{ a}^c.
\end{equation}
Note that the above definition of SCP is very general and does not
include the \textit{standard} requirements; see, for example, \cite
{atar2,atagur,atamanrei,AS}.

\subsection{Static fluid model and throughout sub-optimality}\label{SFM}
The paper deals with certain properties of an underlying fluid model,
to be introduced in this section. We start with the first order
approximations of the parameters.

\begin{assumption} \label{as0}There exist constants $\la_i,\nu_j\in
(0,\infty)$,
$i\in\calI$, $j\in\calJ$ and $\mu_{ij}\in(0,\infty)$,
$(i,j)\in\calE_\mathrm{ a}$, such that
$
n^{-1}\la_i^n\to\la_i, n^{-1}N^n_j\to\nu_j, \mu_{ij}^n\to\mu_{ij}$.
Let $\mu_{ij}=0$ for $(i,j)\in\calE_\mathrm{ a}^c$.
\end{assumption}

The above assumption allows one to imagine a model where arrival and
service processes are
deterministic flows with rates $\la_i$ and
$\mu_{ij}$. There are $J$ service
stations, processing $I$ classes of incoming fluid. Station $j$ has
capacity to hold $\nu_j$
units of fluid. Since routing/scheduling is an important part of
managing the network, an allocation of work among the stations is
pivotal element of the model. The static fluid model uses a fixed
allocation for all times (hence ``static''). Let $\X$ be the set of
\textit{allocation} matrices
\[
\X=\biggl\{ \xi_{ij}, (i,j)\in\calE, \mbox{ such that }
\xi_{ij}\ge0, \mbox { and } \sum_{i\in\calI}
\xi_{ij}\le1, \forall j\in\calJ\biggr\},
\]
where each entry $\xi_{ij}$ represents the fraction
of station's $j$ capacity allocated to process
class-$i$. When station $j$ contains $\psi_{ij}:= \xi_{ij}\nu_j$
units of
class-$i$ fluid, the rate at
which this fluid is processed is $\mu_{ij}\psi_{ij}=\bar\mu_{ij}\xi
_{ij}$, where we set
$\bar\mu_{ij}=\mu_{ij}\nu_j$. The allocation matrix $\xi^*$ to our
model will be chosen according to the following rule.


\begin{assumption}\label{assn1} Consider the following {\emph static
allocation problem \cite{harlop}}:
%
\begin{equation}
\label{25} \min_{\xi\in\X, \rho} \rho, \mbox{ subject to } \sum
_{j\in\calJ}\bar\mu_{ij}\xi_{ij}=
\la_i, \forall i, \sum_{i\in\calI}
\xi_{ij}\le\rho, \forall j,
\end{equation}
and assume it has a unique solution $(\xi^*,\rho^*)$, satisfying:
\begin{longlist}[(1)]
\item[(1)] $\rho^*=1$ and $\sum_{i\in\calI}\xi_{ij}^*=1$ for all
$j\in
\calJ$;
\item[(2)] the set of activities (edges) $(i,j)\in\calE_a$, for which
$\xi^*_{ij}>0$, form a connected tree in a graph with the set of
vertices $\calI\cup\calJ$.
\end{longlist}
\end{assumption}

For convenience, we choose to keep this \textit{standard} set of
assumptions throughout the paper, but, in fact, neither the uniqueness,
nor the tree-like structure is crucial. See more explanation in
Section~\ref{dis}. For the solution $\xi^*$ from Assumption~\ref
{assn1} we denote
%
\begin{equation}
\label{66} \psi_{ij}^*=\xi_{ij}^*\nu_j,\qquad
x^*_i=\sum_j\psi_{ij}^*,\qquad i
\in\calI, j\in\calJ.
\end{equation}
Thus $x^*$
represents the mass of material of each class being processed in all
service stations. The introduced deterministic model, with parameters
$\{ \la, \nu, \mu\}$ and allocation matrix $\{\psi^*\}$, satisfying
Assumption~\ref{assn1}, will be referred to as the \textit{static
fluid model}. Following \cite{harlop,wil}, an activity $(i,j)\in\calE
_\mathrm{ a}$ is said
to be \textit{basic} (resp., \textit{nonbasic}) if $\psi_{ij}^*>0$ (resp.,
$\psi^*_{ij}=0$).

\textit{Throughput sub-optimality.} For $\bar{x}\in\R
^{I}_+$ and $\bar{\nu} \in\R^{J}_+$, define
%
\begin{eqnarray}
\label{XX61}&& \X(\bar x, \bar\nu):= \biggl\{ \psi_{ij}, (i,j) \in
\calE\dvtx \psi _{ij}\ge0,  \sum_{i\in\calI}
\psi_{ij}\le\bar\nu_j, \forall j\in \calJ
\nonumber
\\[-8pt]
\\[-8pt]
\nonumber
&&\hspace*{144pt}\mbox{and }  \sum_{j\in\calJ}\psi_{ij}
\le \bar x_i, \forall i\in\calI \biggr\}.
\end{eqnarray}
Note that from \eqref{66} we have $\psi^*\in\X(x^*, \nu)$. Assumption~\ref{assn1} expresses the
critical load on the system, but
does not discard the possibility that the total processing rate can
exceed the total arrival rate. For a static fluid model we will say
that it is \textit{throughput optimal}
if the following holds:
%
\begin{equation}
\label{50} \mbox{Whenever } \psi\in\X\bigl(x^*,\nu\bigr),
\mbox{one has } \sum
_{(i,j)\in\calE}\mu_{ij}\psi_{ij}\le\sum
_{i\in\calI}\la_i.
\end{equation}

The model is said to be \textit{throughput
sub-optimal} if it is not throughput optimal.

\subsection{The main result}
\label{TO}

The following assumption regards the second order behavior of the
parameters and initial condition.
%

%
\begin{assumption} \label{assn2}There exist $c \in(0,\infty)$,
independent of $n$, such that for $n\ge1$,
%
\begin{eqnarray}\label{1}
\bigl\|n^{-1}\la^n- \la\bigr\| +\bigl\| \mu^n-\mu\bigr\|
+\bigl\|n^{-1} X^n(0) - x^*\bigr\| &\le& cn^{-1/2},
\nonumber
\\[-8pt]
\\[-8pt]
\nonumber
 \bigl\llVert n^{-1} N^n-\nu^n \bigr
\rrVert &\le&(1/2) n^{-1/2}.
\end{eqnarray}
\end{assumption}

\begin{theorem} \label{main_thm}
Let Assumptions \ref{as0}--\ref{assn2} hold. Assume $I=2$ or $J=2$. If,
for some $T>0$, there exists a sequence of SCPs, under which
%
\begin{equation}
\label{weaknc_1} \int_0^T
\bigl\{{e\cdot X^n(s)\ge e\cdot N^n}\bigr\}\,ds \to0\qquad
\mbox{in probability,}
\end{equation}
then the underlying static fluid model is throughput sub-optimal.
\end{theorem}

\section{Proof of Theorem \texorpdfstring{\protect\ref{main_thm}}{2.4}}\label{proof}

\subsection{Intuition and preparations}\label{road} First, we outline
the main ideas of the proof. Fix $n$. It would be convenient to rescale
the system dynamics with respect to the static fluid model. Namely, we
rewrite \eqref{23}--\eqref{40} in the form
%
\begin{eqnarray}
\label{23_sc} \widehat{X}_i^n(t)&=&\widehat{X}_i^{n}(0)+
\widehat{W}_i^n(t) -\sum_{j\in\calJ}
\mu_{ij} \int_0^t\widehat{
\Ps}_{ij}^n(s)\,ds,\qquad i\in\calI,
\\
\label{20_sc} \sum_{j\in\calJ}\widehat
\Ps^n_{ij}(t)&\le&\widehat X_i^n(t),\qquad
i\in \calI; \qquad\sum_{i\in\calI}\widehat\Ps^n_{ij}(t)
\le\widehat N_j^n,\qquad j\in \calJ,
\end{eqnarray}
where we use
%
\begin{eqnarray}\label{scaled}
\nonumber
\widehat{A}^n_i(t)&=&n^{-1/2}
\bigl(A^n_i(t)-\la^n_i t\bigr),\qquad
\widehat {S}^n_{ij}(t)=n^{-1/2}
\bigl(S^n_{ij}(t)-\mu^n_{ij} t\bigr),
\\
  \widehat{X}^n_i(t)&=&n^{-1/2}
\bigl(X^n_i(t)-n x_i^*\bigr),\qquad \widehat{\Psi
}^n_{ij}(t)=n^{-1/2}\bigl(\Psi^n_{ij}(t)-n
\psi^*_{ij}\bigr),
\\
\nonumber
\widehat N^n_j &=& n^{-1/2}
\bigl(N^n_j - n \nu_j \bigr)
\end{eqnarray}
and
%
\begin{eqnarray}
\label{45} \widehat{W}_i^n(t) &=&
\widehat{A}^n_i(t)- \sum_{j\in\calJ}
\widehat {S}_{ij}^n \biggl(\int_0^t
\Ps_{ij}^n(s)\,ds \biggr)
\nonumber
\\[-8pt]
\\[-8pt]
\nonumber
& &{}+n^{-1/2}
\bigl(\la_i^n-n\la_i\bigr)t-n^{-1/2}
\sum_{j\in\calJ}\bigl(\mu _{ij}^n-\mu
_{ij}\bigr)\int_0^t
\Ps_{ij}^n(s)\,ds. 
\end{eqnarray}

The proof will be completed in several steps. The basic principle would
be to show that, once the underlying fluid model is throughput optimal,
it is impossible to quickly eliminate a nonnegligible surplus of customers.

$\bullet$ Our main candidate for a fast unloading of the system will be
the last term of \eqref{23_sc}, since $\widehat{W}^n$ is well known to
be tight; see, for example, \cite{atar2,atamanrei,AS}. Now, due to
throughput optimality \eqref{XX61}, \eqref{50}, since $\Ps^n \in\X
(X^n, N^n )$, we have a crude estimate
%
\begin{equation}
\label{AA9a} \sum_{ij}\mu_{ij} \widehat{
\Ps}^n_{ij}(t) \le\mu_{\mathrm{max}} \bigl(\bigl\|
\widehat{X}^n(t)\bigr\|+\bigl\|\widehat{N}^n\bigr\| \bigr),\qquad t\ge0,
\end{equation}
for $\mu_{\mathrm{max}} = \max_{ij}\{\mu_{ij}\}$, which tells us that, in
principle, the left-hand side of \eqref{AA9a} can be made large by \textit
{quickly} increasing $\| \widehat{X}^n\|$. Of course, stopping the
service (partially or completely) will do the trick, but will not serve
our purpose, thus inviting the question whether, and if so, in what
directions, $ \widehat{X}^n$ can be quickly changed without significant
increase of the total mass $e \cdot\widehat{X}^n$.

$\bullet$ To answer the above we would need Theorem~\ref{12121212} of
Section~\ref{rep}, namely, \textit{representation} \eqref{o12a},
showing that it can be done by using the nonbasic activities along the
so-called zero simple paths, the objects first introduced in \cite{AS},
but with $\mu(p)=0$. To make this paper self contained we have included
Section~\ref{sp}, reminiscing about the basic definitions of simple
paths from \cite{AS} as well as their connection to throughput
optimality (Theorem~\ref{thm}).

$\bullet$ The representation theorem prompts us back to the static
fluid model in an attempt to understand whether one can increase the
throughput by inflicting changes along zero paths. The corresponding
Theorem~\ref{elimin} of Section~\ref{el} provides the \textit{desired}
negative answer and culminates in its dynamic version (Lemma~\ref
{estimate_main} of Section~\ref{key}), \textit{essentially} saying that
there is \textit{no way} to quickly increase $\sum_{ij}\mu_{ij}
\widehat
{\Ps}^n_{ij}$ without increasing $e \cdot\widehat{X}^n$, which is
quite the opposite of what we are trying to achieve.

$\bullet$ The details are finalized in Section~\ref{fin}.

\subsection{Simple paths. Characterization of throughput
sub-optimality}\label{sp}

Denote the index set for all customer classes and service
stations by $\calV:=\calI\cup\calJ$.
For a nonempty set $V$ and $E\subseteq V\times V$, we write $G=(V,E)$ for
the graph with vertex set $V$ and edge set $E$; see, for example, \cite
{dies} for standard definitions. A connected graph that does not
contain cycles is called a $\mathit{tree}$. We denote $\calG_\mathrm{ a}=(\calV,\calE
_\mathrm{ a})$ and refer to it as the graph of activities.

Define the graph of basic activities $\calG_{ ba}$ to be
the subgraph of $\calG_\mathrm{ a}$ having $\calV$ as a vertex set, and
the collection
\[
\calE_{ ba}:=\bigl\{(i,j)\in\calE_\mathrm{ a} \dvtx
\xi^*_{ij}>0\bigr\}
\]
of basic activities as an edge set. By Assumption~\ref{assn1}, the
graph $\calG_{ba}$ is a tree, and by construction of
it as a subgraph of $\calG_a$, all its edges are of the form $(i,j)$
where $i\in\calI$ and $j\in\calJ$. In the definition below and
elsewhere in this section, it will be convenient to identify $(i,j)$
with $(j,i)$ (where $i\in\calI$ and $j\in\calJ$) when referring to
an element of the edge set $\calE$. Although the notation is abused,
there will be no confusion, since $\calI$ and $\calJ$ do not
intersect.

\begin{definition}\label{sim_path}
(i) A subgraph $q=(\calV_q,\calE_q)$ of $\calG_{ba}$ is called a
\textit{basic path} if one has $\calV_q=\{i_0,j_0,\ldots,i_k,j_k\}$ and
\[
\calE_q=\bigl\{(i_0,j_0),(j_0,i_1),
\ldots,(i_k,j_k)\bigr\},
\]
where
$k\ge1$ and $i_0,\dots,i_k\in\calI$, $j_0,\ldots,j_k\in\calJ$ are
$2k+2$ distinct vertices. Note that every edge of a basic path is a
basic activity (i.e., an element of $\calE_{ba}$). Basic paths are used
to define simple paths, as follows:

(ii)
Let the leaves $i_0$ and $j_k$ of a basic path $q$ be denoted by
$i^q$ and, respectively, $j^q$. The pair $(i^q,j^q)$ could be an
activity (an element of $\calE_a$), in which case it is necessarily
a nonbasic activity (i.e., an element of
$\calE_a\setminus\calE_{ba}$), and we say that the graph
$ (\calV_q,\calE_q\cup\{(i^q,j^q)\} )$ is a \textit
{closed simple
path}; otherwise $(i^q,j^q)$ is not an activity (i.e., it is in
$\calE_a^c$), and we say that $q$ itself is an \textit{open simple
path.} We say that $p$ is a \textit{simple path} if it is either a
closed or an open simple path. Let $\mathit{SP}$ be the set of
simple paths.
\end{definition}

\begin{figure}[b]

\includegraphics{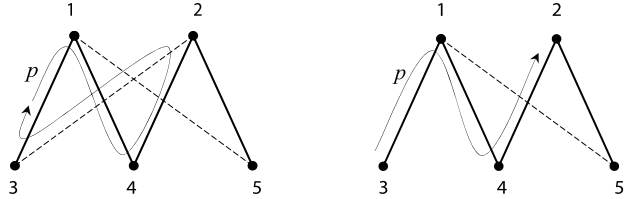}

\caption{Simple paths for cases A and B: On the left $p$
is a closed simple path, while on the right $p$ is open. For case A,
$\mu_{23}>0$ and $(2,3)$ is a
nonbasic activity. For case B, $\mu_{23}=0$ and $(2,3)$ is not an
activity.}\label{fig2}
\end{figure}

\begin{example*} Consider the following static fluid
model, with $2$ classes of customers and $3$
stations:
\[
\nu=\pmatrix{1
\vspace*{2pt}\cr
1 }, \qquad \lambda=\pmatrix{8
\vspace*{2pt}\cr
4 }
\]
and
\begin{eqnarray*}
\mathit{Case\ A:}\quad \mu&=&\bar{\mu}=\pmatrix{ 3 & 10 &
1
\vspace*{2pt}\cr
1 & 4 & 2 };\\
 \mathit{Case\ B:}\quad \mu&=&\bar{\mu}=
\pmatrix{ 3 & 10 & 1
\vspace*{2pt}\cr
0 & 4 & 2}.
\end{eqnarray*}

The resulting optimal static allocation \eqref{25}, \eqref{66} in both
cases is given as
\[
\psi^*=\xi^*=\pmatrix{1 & 0.5 & 0
\vspace*{2pt}\cr
0 & 0.5 & 1 }\quad\mbox{and}\quad x^*=\pmatrix{1.5
\vspace*{2pt}\cr
1.5 },
\]
and we can visualize the graph of activities on Figure~\ref{fig2}. In both cases
we have the same $\calG_{ba}$, consisting of vertices $\{1,2,3,4,5\}$
and edge set $\calE_{ba} =\{(1,3), (1,4), (2,4), (2,5) \}$. Similarly,
both cases have two basic paths [recall, we identify $(i,j)$ with $(j,i)$]
\[
q_1 = \bigl\{(3,1), (1,4), (4,2)\bigr\} \quad\mbox{and}\quad
q_2=\bigl\{(1,4), (4,2), (2,5)\bigr\}.
\]
The basic path $q_1$, together with the corresponding leaves $3$ and
$2$, defines a path $p=\{(3,1), (1,4), (4,2), (2,3)\}$ which will be
closed if $\mu_{23}>0$ [i.e., $(2,3)$ is an activity, case A] and open
otherwise (case B). The only other possible simple path $\{(5,2),
(2,4), (4,1), (1,5)\} $ in both cases will be a closed one.
\end{example*}

Next, we associate directions with edges of simple paths. Let $p$ be
a simple path, and let $q=q^p=(\calV_q,\calE_q)$ be the corresponding
basic path with $\calE_q=\{(i_0,j_0),\ldots,(i_k,j_k)\}$, where
$i_0,\ldots,i_k\in\calI$ and $j_0,\ldots,j_k\in\calJ$. The direction
that will be associated with the edges in $\calE_q$, when considered
as edges of $p$, is as follows: $j_k\to i_k\to j_{k-1}\to
i_{k-1}\to\cdots\to j_0\to i_0$. In the case of an open simple
path, this exhausts all edges of $p$. In the case of a closed simple
path, the direction of $(i_0,j_k)$ is $i_0\to j_k$. We
note that an edge (corresponding to a basic activity) may have
different directions when considered as an edge of different simple
paths. We signify the directions along simple paths by $s(p,i,j)$,
defined for $i\in\calI$, $j\in\calJ$, $(i,j)\in\calE_p$, $p\in \mathit{SP}$,
as
%
\begin{equation}
\label{04} s(p,i,j)=\cases{ -1,&\quad $\mbox{if $(i,j)$, considered as an
edge of
$p$,}$\vspace*{2pt}\cr
&\quad$\mbox{is directed from $i$ to $j$},$\vspace*{2pt}
\cr
  1,&\quad $\mbox{if $(i,j)$,
considered as an edge of $p$,}$\vspace*{2pt}\cr
&\quad $\mbox{is directed from $j$ to $i$.}$}
\end{equation}
Set
%
\begin{equation}
\label{41} m_{i,p} =\sum_{j:(i,j)\in p}s(p,i,j)
\mu_{ij}, \qquad i\in\calI,\qquad  m_p =(m_{i,p},i\in\calI)
\end{equation}
and
%
\begin{equation}
\label{41a} \mu(p)=\sum_{i:(i,j)\in p}m_{i,p}=\sum
_{(i,j)\in\calE
_p}s(p,i,j)\mu _{ij},\qquad  i\in\calI.
\end{equation}

\begin{example*}[(cont.)] Referring to the simple path $p$,
for case A we have $\mu(p) = -7+3=-4$ since
\begin{eqnarray*}
m_{1,p} &=& s(p,1,3)\mu_{13} +s(p,1,4)\mu_{14} =
\mu_{13}-\mu_{14} = -7,
\\
m_{2,p}& =& s(p,2,3)\mu_{23} +s(p,2,4)\mu_{24} =-
\mu_{23}+\mu_{24} = 3.
\end{eqnarray*}
Similarly, for the case B we have $m_{1,p} = -7$, $m_{2,p} = 4$ (since
$\mu_{23}=0$) and $\mu(p) =-3$.
\end{example*}

\begin{theorem}[(Theorem~2, \cite{AS})]\label{thm}
Let Assumptions \ref{assn1} and \ref{assn2} hold. Then the following
statements are equivalent:
\begin{longlist}[(1)]
\item[(1)] the static fluid model is throughput sub-optimal;
\item[(2)] there exists a simple path $p\in \mathit{SP}$ such
that $\mu(p)<0$.
\end{longlist}
\end{theorem}

\begin{example*}[(cont.)] Both cases have a path with $\mu
(p)<0$, hence both are throughput sub-optimal. To see that, for
example, the fluid model in case A is throughput sub-optimal, let
$\beta
> 0$
be sufficiently small, and consider the allocation matrix
\[
\widehat{\xi} = \pmatrix{1-\beta& 0.5+\beta& 0
\vspace*{2pt}\cr
\beta& 0.5-\beta& 1}.
\]
Clearly, we have
$\sum_j\widehat{\xi}_{ij}\nu_j = x_i^*$ for every $i$. However,
$\sum_{(i,j)\in\calE}\widehat{\xi}_{ij}\bar{\mu}_{ij}>\la_1+\la_2$.
\end{example*}

\subsection{Representation}\label{rep}

\begin{theorem}\label{12121212}
Let $X^n$ and $\Ps^n$ satisfy \eqref{23}--\eqref{40}. Then there exist
processes $\Phi^n$, $M^n$ and $ \Upsilon^n$, satisfying:
\begin{longlist}[(1)]
\item[(1)] $\Phi^n(t) \in\X(X^n(t), N^n) $ for $t\ge0$ and
$
\Phi_{ij}^n\equiv0$ for $(i,j) \notin\calE_{ba}
$;
\item[(2)] $
\|\widehat{\Phi}^n(t) \|\le c_F  (\| \widehat{X}^n(t)\| + \|
\widehat N^n \| )
$
for some constant $c_F$, independent of $t, n$;
\item[(3)] $M^n\in\D (\R^{|\mathit{SP}|}  )$, $ \Upsilon^n\in\D
(\R^{|\calI|}  )$, are component-wise nondecreasing, initially zero,
so that the following holds for the scaled processes\setcounter{footnote}{1}\footnote{ We use
$\widehat{\Phi}^n_{ij}:= n^{-1/2} (\Phi^n_{ij} - n\psi
^*_{ij}
)$, $\widehat M^n:=n^{-1/2}M^n$ and $\widehat\Upsilon
^n:=n^{-1/2}\Upsilon^n$.} for $t\ge0$, $i\in. \calI$
%
\begin{eqnarray}
\label{o12} \widehat X_i^n(t)& = &\widehat
X_i^{n}(0)+\widehat W_i^n(t) -
\sum_{j:
(i,j)\in\calE_a }\mu_{ij}\int_0^t
\widehat{\Phi}_{ij}^n(s)\,ds
\nonumber
\\[-8pt]
\\[-8pt]
\nonumber
&&{} + \sum
_{p\in \mathit{SP}}m_{i,p} \widehat M_p^n(t)
+ \widehat\Upsilon_i^n(t).
\end{eqnarray}
\end{longlist}
\end{theorem}

The proof is relegated to the \hyperref[app]{Appendix}.
Together with inequality (2), which is obviously stronger than \eqref
{AA9a}, the theorem indicates that the last two terms of~\eqref{o12}
are the only possible reasons for the abrupt change of $\|\widehat X^n\|
$. The summation term is associated with simple paths, while the last
term corresponds to direct nonwork conservation; see the proof for
more details. The theorem can be viewed as a generalization of Theorem~2.1 from \cite{AMS} where only
closed simple paths (called cycles) were considered.

For a simple path $p\in \mathit{SP}$, we say $p\in\calP_0$, (resp., $p\in
\calP
_-$; or $p\in\calP_+$) if \mbox{$\mu(p)=0$}, [resp., $\mu(p)<0$; or $\mu
(p)>0$ ]. Depending on the subscript sign of $\calP$ the paths will be
called, respectively, \textit{zero}, \textit{negative} or \textit
{positive} paths.

If the static fluid model is throughput optimal (in which case Theorem~\ref{thm} implies $\mathcal{P}_{-} = \varnothing$ ), we rewrite
\eqref
{o12} as
%
\begin{eqnarray}
\label{o12a} \widehat X_i^n(t) = \widehat
X_i^{n}(0)+\widehat W_i^n(t) -
\sum_{j:
(i,j)\in\calE_{ba} }\mu_{ij}\int_0^t
\widehat{\Phi}_{ij}^n(s)\,ds + \widehat\zeta^n_i(t)
+ \widehat\eta^n_i(t),
\nonumber
\\[-8pt]
\\[-8pt]
 \eqntext{i\in\calI,}
\end{eqnarray}
where
%
\begin{equation}
\label{o12b} \quad\widehat\zeta^n_i(t) = \sum
_{p\in\calP_0}m_{i,p} \widehat M_p^n(t),\qquad \widehat\eta^n_i(t)=\sum_{p\in\calP_+}m_{i,p}
\widehat M_p^n(t) +\widehat\Upsilon_i^n(t).
\end{equation}
Notice that $\widehat\zeta^n$ and $\widehat\eta^n$ satisfy [due to
\eqref{41a} and nonnegativity of $\Upsilon_i^n$]
%
\begin{equation}
\label{o12c} e\cdot\widehat\zeta^n (t)\equiv0,\qquad  e\cdot\widehat
\eta^n (t)\ge0.
\end{equation}

\subsection{Discarding zero paths}\label{el}
From \eqref{o12a}--\eqref{o12c} we see that both $\widehat{\zeta}^n$
and $\widehat\eta^n$ can lead to abrupt increase of $\|\widehat X^n\|
$, though only $\widehat\eta^n$ that can do such for $e\cdot\widehat
{X}^n$. The next deterministic (key!) result, viewed as a prelude to
estimate \eqref{p11aa} of Lemma~\ref{estimate_main}, discards any
significant influence of zero paths (represented by $\widehat\zeta^n$)
on system's drift.

\begin{theorem}\label{elimin}
Assume that the static fluid model, (as defined in Section~\ref{SFM}),
is throughput optimal. Take an arbitrary vector $M\in\R_+^{|\calP_0|}$,
with $\|M\|$ small enough,
%
and set
%
\begin{equation}
\label{AA4} x = x^* +\sum_{p\in\mathcal{P}_{0}}m_{p}
M_p.
\end{equation}
Then, if either $I=2$ or $J=2$, the following inequality is true:
%
\begin{equation}
\label{AA5} \sum_{ij}\mu_{ij}
\psi_{ij} \le\sum_{ij}\mu_{ij}
\psi_{ij}^*
\end{equation}
for all $\psi\in\X(x,\nu)$.
\end{theorem}

Before proving the theorem, we point out an important corollary.

\begin{corollary} \label{cor_el} Let the static fluid model be
throughput optimal. Assume we are given some $x_0 \in\R^I_+$,
$\widetilde{\nu} \in\R^J_+$, $\gamma\in\R^I_+$ and a set of numbers
$\{ M_p \ge0, p\in\mathcal{P}_{0} \cup\mathcal{P}_{+}\}$ with $\|M
\|$ sufficiently small. Define $\widetilde{x}= x_0 + \zeta+\eta$,
where
%
\begin{equation}
\label{o12dd} \zeta_i = \sum_{p\in\calP_0}m_{i,p}
M_p,\qquad \eta_i:=\sum_{p\in
\calP
_+}m_{i,p}
M_p+ \gamma_i \qquad \forall i\in\calI.
\end{equation}
Then, if either $I=2$ or $J=2$, for all $\psi\in\X(\widetilde{x},
\widetilde{\nu})$, we have
%
\begin{equation}
\label{AA3} \sum_{ij}\mu_{ij}\bigl(
\psi_{ij} -\psi_{ij}^*\bigr)\le c_{\mu} \bigl( \bigl\|
x_0-x^*\bigr\| + \|\widetilde{\nu}-\nu\| +e\cdot\eta \bigr),
\end{equation}
where $c_{\mu}$ is a constant, independent of $\xi$, $\eta$, $M$.
\end{corollary}

\begin{pf}
Just note that
$
\widetilde{x}=x^* + \zeta+ (x_0-x^*) +\eta= x+(x_0-x^*) +\eta
$
for $x$ from \eqref{AA4}, together with \eqref{AA5}\vspace*{1.5pt} yielding $\sum_{ij}\mu_{ij}(\psi_{ij} -\psi_{ij}^*)\le\mu_{\mathrm{max}}  ( \|
x_0-x^*\|+
\|\widetilde{\nu}-\nu\| +\| \eta\| )$ and the corollary follows
since $\mu(p) > 0$ for each $p\in\calP_+$, and
%
\begin{eqnarray}\label{qqq}
\mu_{\mathrm{max}} \|\eta\| &\le&\mu_{\mathrm{max}}  \biggl(\sum
_{p\in\calP_+}\| m_{p}\| M_p + e\cdot
\gamma_i \biggr)
\nonumber
\\[-8pt]
\\[-8pt]
\nonumber
&\le& c_{\mu}  \biggl(\sum_{p\in\calP_+}
\mu(p) M_p + e\cdot\gamma_i \biggr) = c_{\mu} (e
\cdot \eta)
\end{eqnarray}
for
$
c_{\mu} = \mu_{\mathrm{max}}  (1+ \min \{c\ge0 \dvtx\|m_p\| \le c
\mu(p),
\mbox{for all} p\in\calP_+  \} )$.
\end{pf}

\begin{pf*}{Proof of Theorem~\ref{elimin}} We will start with a basic case
when $I=J=2$ then extend it to more general systems.

\textit{Case} 1: let $\calI= \{1,2\}$ and $\calJ=\{3,4\}$, and assume
the (unique) basic path is given as $q = \{ (3,1), (1,4), (4,2) \}$
with $(2,3)$ being either nonbasic activity or not an activity. The
corresponding simple path $p$ belongs to $\calP_0$, and hence satisfies
%
\begin{equation}
\label{BB1} m_{1,p}+m_{2,p}=(\mu_{13}-
\mu_{14}) + (\mu_{24}-\mu_{23})=0.
\end{equation}
Take a small enough $M>0$, set $\Delta:= (\mu_{13}-\mu_{14})M= (\mu
_{23}-\mu_{24})M$ and define a new $x=(x_1, x_2) = ( x_1^* + \Delta,
x_2^* - \Delta)$. Because of \eqref{BB1}, an elementary argument
implies that any \textit{throughput optimizing} allocation matrix
$\psi
$ is of the form $\psi=\psi^{\gamma}$
%
\begin{equation}
\label{BB2} \bigl(\psi_{13}^{\gamma}, \psi_{14}^{\gamma},
\psi_{23}^{\gamma}, \psi _{24}^{\gamma} \bigr) =
\bigl(\psi_{13}^* -\gamma, \psi_{14}^* +\Delta+ \gamma,
\gamma, \psi_{24}^* -\Delta- \gamma\bigr)
\end{equation}
for some $0\le\gamma\le\min\{\psi^*_{13}, \psi_{24}^*-\Delta\}$,
with the total throughput remaining a constant, independent of $\gamma$,
%
\begin{equation}
\label{BB3} \sum_{ij}\psi_{ij}^{\gamma}
\mu_{ij} \equiv\psi^*_{13} \mu_{13}+ \bigl(\psi
_{14}^*+\Delta\bigr) \mu_{14}+\bigl(\psi_{24}^*-
\Delta\bigr) \mu_{24}.
\end{equation}
%

Assume, on the contrary, that $\sum_{ij}\psi_{ij}^{\gamma} \mu_{ij}
>\la_1 +\la_2$. Due to \eqref{BB3}, the latter inequality will hold for
any feasible $\gamma$. In particular, take $\gamma_0=M \mu_{24}$. It is
easy to check that with such a choice, we have [recall $\Delta= M(\mu
_{23}-\mu_{24})$]
%
\begin{equation}
\label{BB6}  \psi_{23}^{\gamma_0} \mu_{23} +
\psi_{24}^{\gamma_0} \mu_{24} = \la_2.
\end{equation}
Together with $\sum_{ij}\psi_{ij}^{\gamma_0}\mu_{ij} > \la_1 +\la_2$,
it means
%
\begin{equation}
\label{BB7} \psi_{13}^{\gamma_0}\mu_{13} +
\psi_{14} ^{\gamma_0} \mu_{14} >\la_1,
\end{equation}
clearly contradicting the static fluid allocation problem; see
Assumption~\ref{assn1}. Indeed, \eqref{BB6}--\eqref{BB7} means there is
a static fluid allocation $(\widetilde{\psi}_{13}$, $\psi_{14}
^{\gamma
_0}$, $\psi_{23}^{\gamma_0}$, $ \psi_{24}^{\gamma_0} )$, with
$\widetilde{\psi}_{13} < \psi_{13}^{\gamma_0}$, that \textit{fully
serves} each of the two incoming classes without using all the
capacity. 

\textit{Case} 2: now consider the case $I=2$ or $J=2$. An important
property of such systems is that each simple path consists of four
vertices and three or four edges, depending whether or not it is open
or closed; and the arguments from case 1 will be very helpful. In
particular, we argue that the statement of the theorem remains true if
{ \it only one} zero path modification is applied, that is, if $x = x^*
+ m_{p} M_p$ for some path $p\in\mathcal{P}_{0}$, then $\sum_{ij}\mu
_{ij}\psi_{ij} \le\sum_{ij}\mu_{ij}\psi_{ij}^*$ for any $\psi\in
\X
(x,\nu)$. Indeed, let, on the contrary, there exist a throughput
maximizing allocation $\psi\in\X(x,\nu)$ satisfying $\sum_{ij}\mu
_{ij}\psi_{ij} > \sum_{ij}\mu_{ij}\psi_{ij}^*$. Let $\mathcal
{V}_p=\{
i_1, i_2, j_1, j_2 \}$ with a nonbasic $(i_2, j_1)$. Then, again, due
to $\mu_{i_1,j_1}-\mu_{i_1, j_2}+\mu_{i_2,j_2}-\mu_{i_2, j_1}=0$
[recall \eqref{41a} that $p\in\calP_0$], we have that the following
allocation:
\begin{eqnarray}
&&\psi^o_{ij} = \psi_{ij}^*\qquad \mbox{for } (i,j)
\notin p,\nonumber
\\
\label{BB2a}&& \bigl(\psi_{i_1,j_1}^{o}, \psi_{i_1,j_2}^{o},
\psi_{i_2, j_1}^{o}, \psi _{i_2,j_2}^{o} \bigr)
\nonumber
\\[-8pt]
\\[-8pt]
\nonumber
&&\qquad=
\bigl(\psi_{i_1,j_1}^* -\gamma, \psi_{i_1,j_2}^* +\Delta + \gamma,
\gamma, \psi_{i_2,j_2}^* -\Delta- \gamma\bigr)
\end{eqnarray}
will\vspace*{1pt} satisfy $\sum_{ij}\mu_{ij}\psi_{ij} ^o= \sum_{ij}\mu_{ij}\psi
_{ij}> \sum_{ij}\mu_{ij}\psi_{ij}^*$ for any feasible $\gamma$ and
$\Delta= M_p m_{i_1, p} = - M_p m_{i_2, p}$, bringing us precisely to
the first case and, hence, to a contradiction.

Now we extend the latter to several zero paths. Set $k=|\mathcal
{P}_{0}|>1$. Once again, assume that there exists a throughput
maximizing matrix $\psi\in\X(x,\nu)$ that satisfies $\sum_{ij}\mu
_{ij}\psi_{ij} > \sum_{ij}\mu_{ij}\psi_{ij}^*$. Consider an allocation
matrix $\overline{\psi}$ of the form $\overline{\psi}_{ij} = \sum_{p\in
\mathcal{P}_{0}} \psi^{ (p)}_{ij}$, where [slightly abusing the
notation and denoting $\mathcal{V}_p=\{ i_1^p, i_2^p, j_1^p, j_2^p \}$
with a nonbasic $(i_2^p, j_1^p)$ per each path $p$],
%
\begin{eqnarray}
\label{BB2b} && \psi^{ (p)}_{ij} = \frac
{1}{k}
\psi^*_{ij}\qquad \mbox{for } (i,j)\notin\mathcal{E}_{p},
\nonumber\\
&&\bigl( \psi^{(p)}_{i_1^p,j_1^p}, \psi^{ (p)}_{i_1^p,j_2^p},
\psi^{
(p)}_{i_2^p, j_1^p}, \psi^{ (p)}_{i_2^p,j_2^p} \bigr)\\
&&\qquad =
\biggl(\frac
{1}{k}\psi_{i_1^p,j_1^p}^*, \frac{1}{k}
\psi_{i_1^p,j_2^p}^* +\Delta ^p, 0, \frac{1}{k}
\psi_{i_2^p,j_2^p}^* -\Delta^p \biggr),\nonumber
\end{eqnarray}
with $\Delta^p = M_p m_{i_1^p, p} = - M_p m_{i_2^p, p}$. Once again,
since each simple path $p$ belongs to $\calP_0$, we have $\sum_{ij}\mu
_{ij}\overline{\psi}_{ij} = \sum_{ij}\mu_{ij}\psi_{ij}> \sum_{ij}\mu
_{ij}\psi_{ij}^*$. Now consider $k$ \textit{completely separated from
each other} systems with identical set $\{\mu_{ij}\}$, but with arrival
rates and capacities divided by $k$. Clearly, the values $\{\frac
{1}{k}\psi^*_{ij}\}$ will solve the static fluid allocation problem in
the smaller systems. Let each of the smaller systems correspond to each
of the possible $p \in\mathcal{P}_{0}$. To each system apply a
modification along the corresponding path
%
\begin{equation}
\label{BB20} x^{(p)} = \frac{1}{k}x ^* + m_{p}
M_{p}.
\end{equation}
The allocation $\{\psi_{ij}^{ (p)}\}$ from \eqref{BB2b} optimizes the
throughput in the corresponding small system and satisfies (since we
have already treated the case when only one $p\in\calP_0$ has been activated)
%
\begin{equation}
\label{BB21} \sum_{ij}\mu_{ij}
\psi_{ij}^{(p)} \le\frac{1}{k}\sum
_{ij}\mu_{ij} \psi_{ij}^*,
\end{equation}
implying overall
\[
\sum_{ij}\mu_{ij}\overline{
\psi}_{ij} = \sum_{ij}\mu_{ij}
\biggl(\sum_{p\in\mathcal{P}_{0}}\psi_{ij}^{(p)}
\biggr)=\sum_{p\in\mathcal
{P}_{0}} \biggl(\sum
_{ij}\mu_{ij}\psi_{ij}^{(p)}
\biggr)\le\sum_{ij}\mu _{ij}
\psi_{ij}^*,
\]
which completes the proof by contradiction.
\end{pf*}

\subsection{Important estimate}\label{key}

%

Consider the event $\Omega^n_w=\{\| \widehat A^n\|^*_1+\|\widehat S^n\|
^*_1\le5\}$.

\begin{lemma}\label{estimate_main} Let Assumptions \ref{as0}--\ref
{assn2} hold, assume that the static fluid model is throughput optimal,
and let $I=2$ or $J=2$. Then, on the event $\Omega^n_w$, for any
scheduling policy, we have, for $\eps>0$ small enough and $t \le2\eps$,
%
\begin{equation}
\label{p11aa} \sum_{ij}\mu_{ij}\widehat{
\Ps}^n_{ij}(t) \le\eps^{-2/3} \bigl(1+\bigl|\bigl(e\cdot
\widehat{X}^n\bigr)^+\bigr|^*_{t} \bigr).
\end{equation}
\end{lemma}

%
\begin{remark}
In fact, the above inequality holds for some constant $\kappa$, but for
our purposes a crude bound of $\kappa< \eps^{-2/3}$ will be enough as
it saves us the trouble of adjusting {\sl essentially irrelevant}
constants after each operation.
\end{remark}


\begin{pf*}{Proof of Lemma~\ref{estimate_main}}
We will start by showing the relation
%
\begin{equation}
\label{p11} \sum_{ij}\mu_{ij}\widehat{
\Ps}^n_{ij}(t) \le\eps^{-1/2} \bigl(1+\bigl(e\cdot
\widehat{X}^n(t)\bigr)^++ e\cdot\widehat{\eta}^n(t)
\bigr),\qquad t\ge0.
\end{equation}
Recall \eqref{AA9a}
%
\begin{equation}
\label{AA9aa} \sum_{ij}\mu_{ij} \widehat{
\Ps}^n_{ij}(t) \le\mu_{\mathrm{max}} \bigl(\bigl\|
\widehat{X}^n(t)\bigr\|+\bigl\|\widehat{N}^n\bigr\| \bigr),\qquad t\ge0.
\end{equation}
Due to \eqref{1} and since $e\cdot\widehat{\eta}^n$ and is a
nondecreasing process starting at zero, inequality \eqref{p11} will
follow for all $t$ when $\|\widehat{X}^n(t)\|\le\eps^{-1/3}(1+e\cdot
\widehat{\eta}^n(t))$. Now consider the case when $\|\widehat
{X}^n(t)\|
> \eps^{-1/3}(1+e\cdot\widehat{\eta}^n(t))$.

First, assume there is only one class $i$ with $|\widehat{X}_i^n(t)| >
\frac{\eps^{-1/3}}{I}(1+e\cdot\widehat{\eta}^n(t))$. If $\widehat
{X}_i^n(t)<0$, relation \eqref{p11} clearly follows from \eqref{AA9aa}
since the left-hand side of~\eqref{p11} would only increase if
$\widehat
{X}_i^n(t)$ is increased to $-\frac{\eps^{-1/3}}{I}(1+e\cdot\widehat
{\eta}^n(t))$.
Otherwise, if $\widehat{X}_i^n(t)>0$, relation \eqref{p11} follows from
\[
\bigl(e\cdot\widehat{X}^n(t)\bigr)^+\ge e\cdot\widehat{X}^n(t)
\ge\widehat {X}_i^n(t) - \eps^{-1/3}\bigl(1+e
\cdot\widehat{\eta}^n(t)\bigr),
\]
because for all other classes $j\ne i$ we have $\widehat{X}^n_j(t) \ge
-|\widehat{X}^n_j(t)|\ge-\frac{\eps^{-1/3}}{I}(1+e\cdot\widehat
{\eta}^n(t))$.

For the rest of the proof assume that $|\widehat{X}_i^n(t)| > \frac
{\eps^{-1/3}}{I}(1+e\cdot\widehat{\eta}^n(t))$ for several different
$i$'s. From \eqref{o12}--\eqref{o12c} we have (using the fact $t \le
2\eps$)
%
\begin{eqnarray}
\label{AA10} \bigl\| \widehat{X}^n\bigr\|^*_t &\le&\bigl\|
\widehat{X}^n(0)\bigr\|+\bigl\|\widehat{W}^n\bigr\| ^*_{t}+2
\eps c_F \bigl(\bigl\| \widehat{X}^n\bigr\|^*_t +\bigl\|
\widehat{N}^n\bigr\| \bigr)+ \bigl\|\widehat{\zeta}^n(t)\bigr\|
\nonumber
\\[-8pt]
\\[-8pt]
\nonumber
&&{} +\bigl \|
\widehat{\eta}^n(t)\bigr\|.
\end{eqnarray}
Using \eqref{qqq}, we have $\|\eta^n(t)\| \le(c_{\mu}/\mu_{\mathrm{max}})
e\cdot\eta^n(t)$. Moreover, due to the lemma's assumptions, we have
[see \eqref{45}] $\|\widehat{X}^n(0)\| + \|\widehat{W}^n\|^*_{\eps}
\le\eps^{-1/6}$ for $\eps$ small enough, altogether implying
%
\begin{equation}
\label{AA11} \bigl\| \widehat{X}^n(t)\bigr\| \le\bigl\| \widehat{X}^n
\bigr\|^*_t \le\eps^{-1/6} \bigl( 1+\bigl\|\widehat{
\zeta}^n(t)\bigr\|+e\cdot\widehat{\eta}^n(t) \bigr).
\end{equation}
Since $\|\widehat{X}^n(t)\| > \eps^{-1/3}(1+e\cdot\widehat{\eta
}^n(t))$, inequality \eqref{AA11} would imply
%
\begin{equation}
\label{bbbb00} \bigl\|\widehat{\zeta}^n(t)\bigr\| \ge\bigl(
\eps^{-1/6}-1\bigr) \bigl(1+e\cdot\widehat{\eta}^n(t)\bigr),
\end{equation}
that is, there is at least one large ``zero path'' (i.e., $p\in\calP
_0$) activity usage and we are going to apply Corollary~\ref{cor_el} to
``filter out'' the effect of such.

First, if $I>2, J=2$, then all vertices $i \in\mathcal{I}$, except for
one (denote it by $k$), are leaves in the tree of basic activities
$\mathcal{G}_{ba}$. For each leaf
$i_0$ there is a unique simple path $p$, going through $i_0$ and $k$.

Consider the following procedure: Let $\mathcal{I}_0 = \mathcal
{I}_0(t) = \{ i \in\mathcal{I} \setminus\{k\} \dvtx|\widehat
{X}_{i_0}^n(t)| >\frac{ \eps^{-1/3}}{I}(1+e\cdot\widehat{\eta
}^n(t)) \}
$. For $i\in\mathcal{I}_0(t)$ define $\widehat{x}_{i}: =\frac
{\widehat
{X}_{i}^n(t)}{|\widehat{X}_{i}^n(t)|}\frac{ \eps^{-1/3}}{I}(1+e\cdot
\widehat{\eta}^n(t))$, and set $ \widehat{x}_{k}:= \widehat
{X}_{k}^n(t)+\sum_{i\in\mathcal{I}_0(t)} (\widehat{X}_{i}^n(t) -
\widehat{x}_{i}  )$. Finally, for $i \notin(\mathcal{I}_0\cup
\{
k\}) $, set $\widehat{x}_i = \widehat{X}_{i}^n(t)$. Viewing vector
$\widehat{X}^n$ as if it has been obtained from $\widehat{x}$ by
applying $|\mathcal{I}_0|$ zero paths to the latter [as \eqref{1} we
obviously have $\|\widehat{X}^n\|^*_{\eps} \le\|\widehat X^n(0)\|+\|
\widehat{A}^n\|^*_1+2 c n^{1/2}\eps\le\frac{n^{-1/2}}{|\mathit{SP}|} (n\min_{i,j}\psi^*_{ij})$ on $\Omega^n_w$, so the perturbation is indeed
\textit{sufficiently small} when viewed on the fluid level], one can
use Corollary~\ref{cor_el} to get
%
\begin{eqnarray}
\label{CC1} \qquad\sum_{ij}\mu_{ij}\widehat{
\Ps}^n_{ij}(t) &\le&c_{\mu} \bigl( \bigl\| \widehat
{N}^n\bigr\|+\| \widehat{x}\|+e\cdot\widehat{\eta}^n(t) \bigr)
\nonumber
\\[-8pt]
\\[-8pt]
\nonumber
&\le& c_{\mu} \bigl( \bigl\|\widehat{N}^n\bigr\|+
\widehat{x}_{k}^++\eps^{-1/3} \bigl(1+e\cdot\widehat{
\eta}^n(t)\bigr)+e\cdot\widehat{\eta}^n(t) \bigr).
\end{eqnarray}
In the last inequality we once again use the fact that only strictly
positive $\widehat{x}_k$ was worth considering [otherwise the left-hand
side of \eqref{CC1} would only increase if $\widehat{x}_k$ is increased
to $-\frac{\eps^{-1/3}}{I}(1+e\cdot\widehat{\eta}^n(t))$]. A crude
estimate $\widehat{x}_{k}^+ \le (e\cdot\widehat{X}^n(t) )^+
+ \eps^{-1/3}(1+e\cdot\widehat{\eta}^n(t))$ that follows from the
definition of $\widehat{x}_i$ and the relation $\widehat{x}_k= e\cdot
\widehat{X}^n(t)- \sum_{i\in\mathcal{I}_0} \widehat{x}_{i} $
completes the proof of \eqref{p11}. If $I=2$, the same procedure is
applied only once, along any of several possible zero paths. This
proves \eqref{p11}.

To finalize the lemma, note that $\Phi^n$ from \eqref{o12}--\eqref
{o12c} satisfies $\Phi^n(t)\in\X(X^n(t), N^n)$ for all $t$ in the
given range, hence is subject to \eqref{p11} as well. Using that, we have
%
\begin{eqnarray}
\label{bob} \bigl|\bigl(e\cdot\widehat{X}^n\bigr)^+\bigr|^*_{t}
&\ge& e \cdot\widehat{X}^n(0)  - e\cdot \widehat{W}^n (t) +
e \cdot\widehat{\eta}^n(t)
\\
\label{boba} &&{} - \eps^{1/2} \bigl(1+\bigl|\bigl(e\cdot
\widehat{X}^n\bigr)^+\bigr|^*_{t}+e \cdot \widehat {
\eta}^n(t) \bigr)
\\
&\ge& C + \tfrac{1}{2} e \cdot \widehat{\eta}^n(t)-
\eps^{1/2}\bigl|\bigl(e\cdot\widehat{X}^n\bigr)^+
\bigr|^*_{t}
\end{eqnarray}
implying
%
\begin{equation}
\label{bobb} e \cdot\widehat{\eta}^n(t) \le\eps^{-1/5}
\bigl(1+ \bigl|\bigl(e\cdot \widehat {X}^n\bigr)^+\bigr|^*_{t}
\bigr),
\end{equation}
and we complete the proof by substituting \eqref{bobb} into \eqref
{p11}.
\end{pf*}

\subsection{Finalizing the proof}
\label{fin}
%

For arbitrary $\eps>0$, consider the event
%
\begin{eqnarray}\label{omega_3}
\nonumber
\Omega_1^n= \Omega_1^n(
\eps)&=& \Omega_w^n \cap \bigl\{e\cdot \widehat
{X}^n(0)-e\cdot\widehat{N}^n + e\cdot
\widehat{A}^n(\eps)\ge4 \bigr\}
\\
& &{} \cap\bigl\{\bigl \| \widehat{A}^n(\cdot)- \widehat
{A}^n(\eps)\bigr\|^*_{[\eps, 2\eps]} \le1/4 \bigr\}
\\
\nonumber
&&{}\cap \bigl\{ \bigl\| \widehat{S}^n \bigr\|^*_1 \le1/4
\bigr\}.
\end{eqnarray}
It is standard (e.g., Theorem~14.6 in \cite{bil}) that component-wise
both $\widehat A^n$ and $\widehat S^n$ converge weakly to \textit
{independent} Brownian motion processes. Therefore there exist
constants $n_1=n_1(\eps) \in\N$ and $\delta=\delta(\eps)>0$, so that
$\PP (\Omega_1^n )>\delta$ for all $n\ge n_1$.

Fix $\eps>0$. Theorem~\ref{main_thm} guarantees that there exists a
sequence of SCPs satisfying
%
\begin{equation}
\label{E4} \lim_{n\to\infty}\PP \biggl(\Omega_1^n
\cap \biggl\{\int_0^T 1 \bigl
\{{e\cdot X^n(s)\ge e\cdot N^n}\bigr\}\,ds >\eps \biggr\}
\biggr)=0.
\end{equation}
%
Let $\Omega^n = \Omega_1^n \cap \{\int_0^T
1 \{{e\cdot X^n(s)\ge e\cdot N^n}\}\,ds\le\eps \}
$. Relation \eqref{E4} implies that there exists a constant $n_0(\eps
)\in\N$ so that
%
\begin{equation}
\label{E5} \PP \bigl(\Omega^n \bigr)>\frac{\delta}{2}\qquad \mbox{for all $n\ge n_0$}.
\end{equation}
In what follows we assume that the static fluid model is \textit
{throughput optimal} and will come to a conclusion that the event
$\Omega^n $ is impossible (i.e., $\Omega^n$ is an empty set) for
$n\ge
n_0$ and $\eps$ small enough, thus contradicting \eqref{E5}.

From \eqref{23_sc}, \eqref{45}, \eqref{1}, Lemma~\ref
{estimate_main} and \eqref{omega_3} on the event $\Omega^n$,
%
\begin{equation}
\label{xx6} e\cdot\widehat{X}^n(\eps)-e\cdot\widehat{N}^n
\ge7/2- c\eps- \eps ^{1/3} \bigl(1+\bigl|\bigl(e\cdot\widehat{X}^n
\bigr)^+\bigr|^*_{\eps} \bigr),
\end{equation}
which, for $\eps$ small enough, yields
%
\begin{equation}
\label{xx6c} e\cdot\widehat{X}^n(\eps)-e\cdot\widehat{N}^n+
\eps^{1/3}\bigl|\bigl(e\cdot \widehat{X}^n\bigr)^+\bigr|^*_{\eps}
\ge 2,
\end{equation}
giving us two possible scenarios: $e\cdot\widehat{X}^n(\eps)-e\cdot
\widehat{N}^n \ge\eps^{1/3}|(e\cdot\widehat{X}^n)^+|^*_{\eps} $ and
$\eps^{1/3}|(e\cdot\widehat{X}^n)^+|^*_{\eps} \ge e\cdot\widehat
{X}^n(\eps)-e\cdot\widehat{N}^n$.

\textit{Case} 1. Assume $e\cdot\widehat{X}^n(\eps)-e\cdot\widehat
{N}^n \ge\eps^{1/3}|(e\cdot\widehat{X}^n)^+|^*_{\eps} $.\vspace*{1pt} Together
with \eqref{xx6c}, this implies $e\cdot\widehat{X}^n(\eps)-e\cdot
\widehat{N}^n\ge1$. Let $\tau_{\eps} = \inf\{ t>\eps\dvtx e\cdot
\widehat
{X}^n(t)= e\cdot\widehat{N}^n\}$. Notice that $\tau_{\eps}$ is well
defined since the jumps of $e\cdot X^n$ are of size $1$ and, moreover,
satisfies $\tau_{\eps}<2\eps$ on $\Omega^n$, because the total queueing
time does not exceed $\eps$. Using $e\cdot\widehat{X}^n(\eps
)-e\cdot
\widehat{N}^n\ge1$, \eqref{23_sc}, Lemma~\ref{estimate_main} and
\eqref
{omega_3} we can write
\begin{eqnarray*}
0&= & e\cdot\widehat{X}^n(\tau_{\eps})-e\cdot
\widehat{N}^n
\\
& \ge& e\cdot\widehat{X}^n(\eps) -e\cdot
\widehat{N}^n +e \cdot\widehat{W}^n(\tau_{\eps}
\dvtx\eps) -\eps ^{1/3} \bigl(1+\bigl|\bigl(e\cdot\widehat{X}^n
\bigr)^+\bigr|^*_{\tau_{\eps}} \bigr)
\\
&\ge&1/8- \eps^{1/3}\bigl|\bigl(e\cdot\widehat{X}^n
\bigr)^+\bigr|^*_{\tau_{\eps}},
\end{eqnarray*}
implying
%
\begin{equation}
\label{xx4} \bigl|\bigl(e\cdot\widehat{X}^n\bigr)^+\bigr|^*_{\tau_{\eps}}
\ge\eps^{-1/4}.
\end{equation}
In other words, a large queue of at least $ \eps^{-1/4}$ has to be
eliminated before time $\tau_{\eps}$. Let $\alpha$ be the \textit{last}
time before $\tau_\eps$, satisfying $|(e\cdot\widehat
{X}^n)^+|^*_{\tau
_{\eps}}=e\cdot\widehat{X}^n(\alpha)\ge\eps^{-1/4}$. We have
\begin{eqnarray*}
0&= & e\cdot\widehat{X}^n(\tau_{\eps})-e\cdot
\widehat{N}^n
\\
& \ge& e\cdot\widehat{X}^n(\alpha) -e\cdot
\widehat{N}^n +e \cdot\widehat{W}^n(\tau_{\eps}
\dvtx\alpha) -\eps ^{1/3} \bigl(1+e\cdot\widehat{X}^n(\alpha)
\bigr)
\\
&\ge& C+ \tfrac{1}{2} e\cdot\widehat{X}^n(\alpha) \ge C
+ (1/2) \eps^{-1/4},
\end{eqnarray*}
for some constant $C$, which is an obvious contradiction for $\eps$
small enough.

\textit{Case} 2. If $\eps^{1/3}|(e\cdot\widehat{X}^n)^+|^*_{\eps}
\ge
e\cdot\widehat{X}^n(\eps)-e\cdot\widehat{N}^n$, then $|(e\cdot
\widehat{X}^n)^+|^*_{\eps} \ge\eps^{-1/3}$ by~\eqref{xx6c}, and the
same considerations as in the previous case can be applied. Let $\alpha
$ be the \textit{last} time before $\eps$, satisfying $|(e\cdot
\widehat
{X}^n)^+|^*_{\eps}=e\cdot\widehat{X}^n(\alpha)\ge\eps^{-1/3}$,
and define
$\tau_{\alpha} = \inf\{ t>\alpha\dvtx e\cdot\widehat{X}^n(t)=
e\cdot
\widehat{N}^n\}$. Then
\begin{eqnarray*}
0&= & e\cdot\widehat{X}^n(\tau_{\alpha})-e\cdot
\widehat{N}^n
\\
& \ge& e\cdot\widehat{X}^n(\alpha) -e\cdot
\widehat{N}^n +e \cdot\widehat{W}^n(\tau_{\alpha}
\dvtx\alpha) -2\eps ^{1/3} \bigl(1+e\cdot\widehat{X}^n(
\alpha) \bigr)
\\
&\ge& C+ \tfrac{1}{3} e\cdot\widehat{X}^n(\alpha) \ge
C + (1/3) \eps^{-1/3},
\end{eqnarray*}
for some constant $C$, giving the contradiction once again. This
concludes the proof of Theorem~\ref{main_thm}. 

\section{General structures}\label{gs} Theorem~\ref{main_thm} shows
that null-controllability is impossible if the underlying fluid model
is throughput optimal. The result is valid for the case $\min\{I,J\}
=2$, and the assumption is crucial for both Theorem~\ref{elimin} and
Lemma~\ref{estimate_main}. How can Theorem~\ref{main_thm} be extended
for general $I$ and $J$, especially, since it is relatively easy to
\textit{numerically check} conditions \eqref{AA4}--\eqref{AA5}
(enough to check separately for each zero path)? We give a partial answer.

\begin{definition}\label{path_dep}
A path $p\in \mathit{SP}$ is called \textsf{class-dependent} if \eqref{04}--\eqref{41a}
%
\begin{equation}
\label{bbbbb41} \sum_{j:(i,j)\in\calE_p}s(p,i,j)\mu_{ij}
=0,\qquad i\in\calI.
\end{equation}
\end{definition}

There are only two summands for each given $i$ in \eqref{bbbbb41}.
Basically, the definition says that for each $i\in\calI$, belonging to
$p$, and two (just these two!) adjacent activities $(i,j_1)$ and $(i,
j_2)$ from the very same path $p$, we must have $\mu_{i,j_1}=\mu
_{i,j_2}$. Similarly, we have the following:

\begin{definition}\label{path_dep2}
A path $p\in \mathit{SP}$ is called \textsf{pool-dependent} if
%
\begin{equation}
\label{bbbb41} \sum_{i:(i,j)\in\calE_p}s(p,i,j)\mu_{ij}
=0,\qquad j\in\calJ.
\end{equation}
\end{definition}

From \eqref{04}--\eqref{41a}, each of the above two types must be a
zero path, that is, $\mu(p)=0$.

\begin{theorem} \label{main_thm_gen}
Let Assumptions \ref{as0}--\ref{assn2} hold, and let $I, J\ge1$.
Assume that the fluid model is throughput optimal and satisfies one of
the following:
\begin{longlist}[(1)]
\item[(1)] has no zero paths, that is, $\calP_0=\varnothing$;
\item[(2)] each $p\in\calP_0$ is either class- or pool-dependent; or,
for small $\kappa>0$,
%
\begin{equation}
\label{bbbb0} \sum_{ij}\mu_{ij}\bigl(
\psi_{ij} -\psi_{ij}^*\bigr)< 0 \qquad\mbox{whenever } \psi\in\X
\bigl(x^* +m_{p} \kappa, \nu\bigr).
\end{equation}
\end{longlist}
Then it is impossible to find $T>0$ and a sequence of SCPs, satisfying
\eqref{weaknc_1}; that is, (weak) null controllability is impossible.
\end{theorem}

\begin{remark}\label{mark} Currently this is as close as we can get
to the conclusion that, in the general case, \eqref{AA4}--\eqref{AA5}
prescinds null controllability (for throughput optimal fluid models).
Apparently, more work is required when \eqref{bbbb0} results in
equality, with path being neither \textsf{class-} nor \textsf{pool-dependent}. We feel, however, that such situations are very rare,
maybe even impossible (and may as well contradict to uniqueness of the
underlying fluid model; see Assumption~\ref{assn1}).
\end{remark}

\begin{remark} Theorem~\ref{main_thm_gen} trivially implies that
null-controllability is also impossible for either one of the following
types of the fluid model:
\begin{longlist}[(1)]
\item[(1)] the service rates depend only on the class type (\textit
{class-dependent}),
%
\begin{equation}
\label{bbbb33} \mu_{ij}=\mu_i,\qquad i\in\calI, j\in\calJ;
\end{equation}
\item[(2)] the service rates depend only on the station type (\textit
{pool-dependent}),
%
\begin{equation}
\label{bbbb32} \mu_{ij}=\mu_j,\qquad i\in\calI, j\in\calJ.
\end{equation}
\end{longlist}
Indeed, in both cases the fluid model is throughput optimal, and all
paths are either class- or pool-dependent.
\end{remark}

\begin{pf*}{Proof of Theorem~\ref{main_thm_gen}} It will be enough to show
that relation \eqref{p11} of Lemma~\ref{estimate_main} remains intact,
as no other part of the proof of Theorem~\ref{main_thm} has any
structure constraints.

\textit{Case} 1. Relation \eqref{p11} trivially follows from the
current proof of\break Lemma~\ref{estimate_main}.

%
\textit{Case} 2. The argument goes exactly as in the proof of Lemma~\ref
{estimate_main}, until it gets to \eqref{bbbb00}, stating that at least
one ``large'' zero path has been activated. Some extra work has to be
done at this point.

\begin{longlist}[A.]
\item[A.] \textit{All zero paths are class-dependent}: in such
case [see \eqref{04}--\eqref{41a}] we have $\|m_p\|=0$ for each zero
path; hence none of these paths has \textit{any} effect on the system.
Once again, \eqref{p11} follows trivially.

\item[B.] \textit{All zero paths are pool-dependent}: assume
that all zero paths satisfy \eqref{bbbb41}. We start with the case when
there is only one zero path $p$. The following estimate will be useful.
From the exact structure of $\Phi^n$ from Theorem~\ref{12121212}, and
\eqref
{o12a}, we have a crude estimate
%
\begin{equation}
\label{popp1} \bigl|\widehat{X}^n_i(t)\bigr | \le
\eps^{-1/3} \bigl(1+\bigl(e\cdot\widehat {X}^n(t)\bigr)^++e
\cdot\widehat{\eta}^n(t) \bigr),\qquad  i\notin\calV_p,
\end{equation}
since no zero paths have been applied to such classes $i$. Now, for
any feasible allocation $ \Psi^n(t)$ consider a unique, \textit
{standard}, \cite{atar1,atar2} allocation $ \phi^n(t) \in\X
(X^n(t), N^n)$ that is zero for nonbasic activities and is work conserving:
$
\min \{(e\cdot\widehat{X}^n(t)-e\cdot\widehat{N}^n), (e\cdot
\widehat{N}^n - \sum_{ij}\widehat{\phi}_{ij}^n(t))  \}=0$.
By throughput optimality of the fluid model [i.e., $\calP
_-=\varnothing
$; see the definitions before \eqref{o12a}], we must have $\sum_{ij}\mu
_{ij}\widehat{\Ps}^n_{ij}(t) \le\sum_{ij}\mu_{ij}\widehat{\phi
}^n_{ij}(t) $. Using \eqref{popp1} and the structure of $\phi^n$,
%
\begin{equation}
\label{pop16} \bigl|\widehat{\phi}_{ij}^n(t)\bigr| \le
\eps^{-1/3} \bigl(1+\bigl(e\cdot\widehat {X}^n(t)\bigr)^++e
\cdot\widehat{\eta}^n(t) \bigr),\qquad (i,j)\notin\calE_p,
\end{equation}
as well as
%
\begin{equation}
\label{pop15} \bigl|\widehat{\phi}_{i_0,j_0}^n(t)\bigr| \le
\eps^{-1/3} \bigl(1+\bigl(e\cdot \widehat {X}^n(t)\bigr)^++e
\cdot\widehat{\eta}^n(t) \bigr)
\end{equation}
for the leaf $(i_0,j_0)$ of the basic simple path, corresponding to
$p$. Together, \eqref{pop16}, \eqref{pop15} imply that for each station
$j\in\calJ\cap\calV_p$, connecting exactly two path edges
$(i_1,j)\in\calE_p$ and $(i_2,j)\in\calE_p$ [this excludes the leaf
from \eqref{pop15}], we also have
%
\begin{equation}
\label{pop17} \widehat{\phi}_{i_1,j}^n(t) + \widehat{
\phi}_{i_2,j}^n(t) \le\eps ^{-1/3} \bigl(1+\bigl(e
\cdot\widehat{X}^n(t)\bigr)^++e\cdot\widehat{\eta}^n(t)
\bigr)
\end{equation}
and due to pool-dependence along the path \eqref{bbbb41} (otherwise it
will not hold!), we get
%
\begin{eqnarray}\qquad
\label{pop18} \mu_{i_1,j} \widehat{\phi}_{i_1,j}^n(t)
+ \mu_{i_2,j} \widehat {\phi }_{i_2,j}^n(t) &= &
\mu_{i_1,j} \bigl( \widehat{\phi}_{i_1,j}^n(t) +
\widehat{\phi}_{i_2,j}^n(t) \bigr)
\nonumber
\\[-8pt]
\\[-8pt]
\nonumber
&\le& \mu_{i_1,j} \eps ^{-1/3} \bigl(1+\bigl(e\cdot
\widehat{X}^n(t)\bigr)^++e\cdot\widehat{\eta}^n(t) \bigr).
\end{eqnarray}
This proves \eqref{p11}. The extension to several pool-dependent paths
is straight\-forward---the only difference being the inclusion of all
pool-dependent activities (possibly more than two), connected by
station $j$, into the left-hand side of \eqref{pop18}. The right-hand side of \eqref{pop18}
will remain the same.

\item[C.] \textit{All zero paths satisfy} \eqref{bbbb0}: by
linearity, there exists a constant $c>0$, such that for any feasible
set $\{M_p\ge0, p\in\calP_0\}$ and $\psi\in\X(x^* + \sum_{p\in
\calP
_0}m_{p} M_p, \nu)$
%
\begin{equation}
\label{pop3} \sum_{ij}\mu_{ij}\bigl(
\psi_{ij} -\psi_{ij}^*\bigr) \le- c \sum
_{p\in
\calP_0} M_p \le- \eps^{1/2} \sum
_{p\in\calP_0}\|m_p\| M_p.
\end{equation}
Note that the first inequality in \eqref{pop3} becomes an equality if
and only if each of the $M_p$ is zero. Applying \eqref{pop3} to
processes from \eqref{o12a} and using Theorem~\ref{12121212}(2), we get
for any feasible allocation $\Ps^n(t)$,
%
\begin{eqnarray}\label{pop13}
\sum_{ij}\mu_{ij}\widehat{
\Ps}^n_{ij}(t)&\le& - \eps^{1/2} \bigl\| \widehat{\zeta
}^n(t)\bigr\|+ \eps c_F \bigl(\bigl\| \widehat{X}^n
\bigr\|^*_t+\bigl\|\widehat{N}^n\bigr\| \bigr)
\nonumber
\\[-8pt]
\\[-8pt]
\nonumber
&&{}+\eps^{-1/3}\bigl(1+e\cdot\widehat{\eta}^n(t)
\bigr).
\end{eqnarray}
Using relation \eqref{AA11}, we continue
%
\begin{eqnarray}
\label{pop20} \sum_{ij}\mu_{ij}\widehat{
\Ps}^n_{ij}(t)&\le&- \eps^{1/2} \bigl\| \widehat {
\zeta}^n(t)\bigr\|+ \eps^{5/6} \bigl( 1 +\bigl\|\widehat{
\zeta}^n(t)\bigr\| +e\cdot \widehat{\eta}^n(t) \bigr)
\nonumber
\\[-8pt]
\\[-8pt]
\nonumber
&&{}+\eps^{-1/3}\bigl(1+e\cdot\widehat {\eta }^n(t)
\bigr)\le
\eps^{-1/3}\bigl(1+e\cdot
\widehat{\eta}^n(t)\bigr),
\end{eqnarray}
where we used $\|\widehat{\zeta}^n(t)\| (-\eps^{1/2}+\eps^{5/6}+\eps^{5/6}
\eps
^{1/6}  ) \le0$ [again, the strict inequality in \eqref{bbbb0} is
crucial for the existence of the ``$-\eps^{1/2}$'' term!]. And \eqref
{p11} follows.

\item[D.] \textit{Finalizing}: we use B and
C to complete the theorem. In particular, once again, the (\textit
{work-conserving, no nonbasic activities} ) allocation $\phi^n$ will
be introduced. After that, the sum $ \sum_{ij}\mu_{ij}\widehat{\phi
}^n_{ij}(t) $ will be decomposed into two different sums: one will
contain all the terms satisfying \eqref{pop16}--\eqref{pop18} [this
will also include the possible intersections of pool-dependent paths
and paths, satisfying \eqref{bbbb0}]; another part of the summation
will satisfy \eqref{pop20}.
This completes the proof.\quad\qed
\end{longlist}\noqed
\end{pf*}


\section{Final remarks}\label{dis}
The text is an attempt to understand whether a \textit{given} static
fluid model is throughput optimal; and some words need to be said
regarding Assumption~\ref{assn1}, in particular, (a) the treelike
structure and (b) the uniqueness of the solution to \eqref{25} in general.

(a) To start with, null controllability is clearly impossible if the
solution to~\eqref{25} does not contain a basic path of the length at
least $3$ (we need to have at least two stations and two classes to be
connected together by basic activities), so some kind of connectivity
should be assumed. When a connected component contains cycles fully
composed by basic activities, one may have trouble defining
weights/directions along simple paths (as was done in Section~\ref
{sp}), yet one thing will remain true: a throughput optimal model can
only have cycles with weight zero. Otherwise, a positive path can
become negative if applied in the other direction, and vice versa. This
brings us back to the very same zero paths, main ingredients of the
current paper.

(b) A mass vector $x^*$, coming from the solution of \eqref{25} is a
key (!) element due to Assumption~\ref{assn2} about the initial
condition. This invites a reasonable question. What if there is another
optimal solution, with the same vector $x^{*}$, but a different graph
structure? This is clearly feasible, although both (or, infinite
number, in that case) possible static fluid models would still be
either all throughput optimal, or sub-optimal altogether, since
definition \eqref{50} does not require any special graph structure.

Now, what if the other solution has a \textit{different} mass vector,
say, $x^{**}$. Is it possible that $x^*$-solution is throughput
optimal, while $x^{**}$ is not? We claim it is not feasible, at least
in the case $I=2$ or $J=2$, with arguments similar to ones in the proof
of Theorem~\ref{elimin} (since $e\cdot x^* = e\cdot x^{**}$). The more
general structure is still to be resolved$\ldots.$

\begin{appendix}\label{app}
\section*{Appendix: Sketch of the proof of Theorem~\texorpdfstring{\protect\ref{12121212}}{3.3}}

Using the scaling
$\widehat{f}=n^{-1/2}f$, introduce auxiliary processes $\widehat{Y}_i^n(t)$,
representing the scaled number of class-$i$ customers that are in the queue
(and not being served) at time $t$, and $\widehat{Z}_j^n(t)$---the scaled
number of servers at station $j$ that are idle at time $t$. Clearly,
we have the following relations:
%
\begin{eqnarray}
\label{20q} \widehat Y_i^n(t)+\sum
_{j\in\calJ}\widehat\Ps^n_{ij}(t)&=&\widehat
X_i^n(t), \qquad i\in\calI,
\\[-2pt]
\label{21q} \widehat Z_j^n(t)+\sum
_{i\in\calI}\widehat\Ps^n_{ij}(t)&=&\widehat
N_j^n, \qquad j\in\calJ.
\end{eqnarray}
The proof can be viewed as generalization of Theorem~2.1 from Atar,
Mandelbaum and Shaikhet \cite{AMS}, whose decomposition used only
closed simple paths (called cycles). In particular the set $\{\widehat
\Psi^n_{ij}, (i,j) \in\calE_a \}$ was decomposed into basic and
nonbasic activities, respectively, $\{\widehat\Psi^n_{ij}, (i,j) \in
\calE_{ba} \}$ and $\{\widehat\Psi^n_{ij}, (i,j) \in\calE_{ba}^c
\}
$, turning \eqref{23_sc} into (see Section~2.3 in \cite{AMS})
\begin{eqnarray*}
\label{x3x1a}
\widehat X_i^n(t) &=& \widehat
X_i^{n}(0)+\widehat W_i^n(t) -
\sum_{j\in\calJ} \mu_{ij} \int
_0^t G_{ij}\bigl(\widehat
X^n(s)-\widehat Y^n(s),\widehat N^n-\widehat
Z^n(s)\bigr)\,ds
\\[-2pt]
&&{}+\sum_{p: p \mathrm{\mbox{-}closed\ simple\ path}} m_{i,p}\int
_0^t\widehat\Ps_p^n(s)\,ds,
\end{eqnarray*}
where $\widehat\Ps_p^n$ corresponds to a unique nonbasic activity,
associated with simple path $p$ and the function $G$, introduced in
\cite{atar1}.

For our purposes, however, that is not enough, since we want to single
out all the terms that can cause an abrupt change of $\widehat X^n$,
and nonwork conservation is exactly what we are looking for, since, a
priori we do not have the relation \mbox{$e \cdot\widehat Y^n \wedge e
\cdot\widehat Z^n =0$}.\vadjust{\goodbreak}

The rest of the proof follows the lines of Theorem~2.1 in \cite{AMS},
with an additional requirement to cover direct nonwork conservation,
that is, situations when $\widehat Y^n_i \wedge\widehat Z^n_j>0$
while $\mu_{ij}>0$, as well as the open simple paths, corresponding to
what we call an indirect nonwork conservation, that is, situations
when $\widehat Y^n_i \wedge\widehat Z^n_j>0$ while $\mu_{ij}=0$. We
leave out the details. 
\end{appendix}

\section*{Acknowledgements} The author wants to thank the anonymous
referees of~\cite{AS} for suggesting the research topic. The
interpretation (or misinterpretation) of their question, resulting in
the current work, is all mine.

%



\printaddresses


\begin{thebibliography}{13}

\bibitem{atar1}
%
\begin{barticle}[mr]
\bauthor{\bsnm{Atar},~\bfnm{Rami}\binits{R.}}
(\byear{2005}).
\btitle{A diffusion model of scheduling control in queueing systems
with many servers}.
\bjournal{Ann. Appl. Probab.}
\bvolume{15}
\bpages{820--852}.
\bid{doi={10.1214/105051604000000963}, issn={1050-5164}, mr={2114991}}
\end{barticle}
%
\bptok{imsref}%
\endbibitem

\bibitem{atar2}
%
\begin{barticle}[mr]
\bauthor{\bsnm{Atar},~\bfnm{Rami}\binits{R.}}
(\byear{2005}).
\btitle{Scheduling control for queueing systems with many servers:
Asymptotic optimality in heavy traffic}.
\bjournal{Ann. Appl. Probab.}
\bvolume{15}
\bpages{2606--2650}.
\bid{doi={10.1214/105051605000000601}, issn={1050-5164}, mr={2187306}}
\end{barticle}
%
\bptok{imsref}%
\endbibitem

\bibitem{atagur}
%
\begin{barticle}[mr]
\bauthor{\bsnm{Atar},~\bfnm{Rami}\binits{R.}} \AND
\bauthor{\bsnm{Gurvich},~\bfnm{Itai}\binits{I.}}
(\byear{2014}).
\btitle{Scheduling parallel servers in the nondegenerate slowdown
diffusion regime: {A}symptotic optimality results}.
\bjournal{Ann. Appl. Probab.}
\bvolume{24}
\bpages{760--810}.
\bid{doi={10.1214/13-AAP935}, issn={1050-5164}, mr={3178497}}
\bptnote{check year}%
\end{barticle}
%
\bptok{imsref}%
\endbibitem

\bibitem{atamanrei}
%
\begin{barticle}[mr]
\bauthor{\bsnm{Atar},~\bfnm{Rami}\binits{R.}},
\bauthor{\bsnm{Mandelbaum},~\bfnm{Avi}\binits{A.}} \AND
\bauthor{\bsnm{Reiman},~\bfnm{Martin~I.}\binits{M.~I.}}
(\byear{2004}).
\btitle{Scheduling a multi class queue with many exponential servers:
Asymptotic optimality in heavy traffic}.
\bjournal{Ann. Appl. Probab.}
\bvolume{14}
\bpages{1084--1134}.
\bid{doi={10.1214/105051604000000233}, issn={1050-5164}, mr={2071417}}
\end{barticle}
%
\bptok{imsref}%
\endbibitem

\bibitem{AMS}
%
\begin{barticle}[mr]
\bauthor{\bsnm{Atar},~\bfnm{Rami}\binits{R.}},
\bauthor{\bsnm{Mandelbaum},~\bfnm{Avi}\binits{A.}} \AND
\bauthor{\bsnm{Shaikhet},~\bfnm{Gennady}\binits{G.}}
(\byear{2006}).
\btitle{Queueing systems with many servers: Null controllability in
heavy traffic}.
\bjournal{Ann. Appl. Probab.}
\bvolume{16}
\bpages{1764--1804}.
\bid{doi={10.1214/105051606000000358}, issn={1050-5164}, mr={2288704}}
\end{barticle}
%
\bptok{imsref}%
\endbibitem

\bibitem{AS}
%
\begin{barticle}[mr]
\bauthor{\bsnm{Atar},~\bfnm{Rami}\binits{R.}} \AND
\bauthor{\bsnm{Shaikhet},~\bfnm{Gennady}\binits{G.}}
(\byear{2009}).
\btitle{Critically loaded queueing models that are throughput suboptimal}.
\bjournal{Ann. Appl. Probab.}
\bvolume{19}
\bpages{521--555}.
\bid{doi={10.1214/08-AAP551}, issn={1050-5164}, mr={2521878}}
\end{barticle}
%
\bptok{imsref}%
\endbibitem

\bibitem{bil}
%
\begin{bbook}[mr]
\bauthor{\bsnm{Billingsley},~\bfnm{Patrick}\binits{P.}}
(\byear{1999}).
\btitle{Convergence of Probability Measures},
\bedition{2nd} ed.
\bpublisher{Wiley},
\blocation{New York}.
\bid{doi={10.1002/9780470316962}, mr={1700749}}
\end{bbook}
%
\bptok{imsref}%
\endbibitem

\bibitem{dies}
%
\begin{bbook}[mr]
\bauthor{\bsnm{Diestel},~\bfnm{Reinhard}\binits{R.}}
(\byear{2000}).
\btitle{Graph Theory},
\bedition{2nd} ed.
\bseries{Graduate Texts in Mathematics}
\bvolume{173}.
\bpublisher{Springer},
\blocation{New York}.
\bid{doi={10.1007/b100033}, mr={1743598}}
\end{bbook}
%
\bptok{imsref}%
\endbibitem

\bibitem{halwhi}
%
\begin{barticle}[mr]
\bauthor{\bsnm{Halfin},~\bfnm{Shlomo}\binits{S.}} \AND
\bauthor{\bsnm{Whitt},~\bfnm{Ward}\binits{W.}}
(\byear{1981}).
\btitle{Heavy-traffic limits for queues with many exponential servers}.
\bjournal{Oper. Res.}
\bvolume{29}
\bpages{567--588}.
\bid{doi={10.1287/opre.29.3.567}, issn={0030-364X}, mr={0629195}}
\end{barticle}
%
\bptok{imsref}%
\endbibitem

\bibitem{harlop}
%
\begin{barticle}[mr]
\bauthor{\bsnm{Harrison},~\bfnm{J.~Michael}\binits{J.~M.}} \AND
\bauthor{\bsnm{L{\'o}pez},~\bfnm{Marcel~J.}\binits{M.~J.}}
(\byear{1999}).
\btitle{Heavy traffic resource pooling in parallel-server systems}.
\bjournal{Queueing Syst.}
\bvolume{33}
\bpages{339--368}.
\bid{doi={10.1023/A:1019188531950}, issn={0257-0130}, mr={1742575}}
\end{barticle}
%
\bptok{imsref}%
\endbibitem

\bibitem{mansto}
%
\begin{barticle}[mr]
\bauthor{\bsnm{Mandelbaum},~\bfnm{Avishai}\binits{A.}} \AND
\bauthor{\bsnm{Stolyar},~\bfnm{Alexander~L.}\binits{A.~L.}}
(\byear{2004}).
\btitle{Scheduling flexible servers with convex delay costs:
Heavy-traffic optimality of the generalized {$c\mu$}-rule}.
\bjournal{Oper. Res.}
\bvolume{52}
\bpages{836--855}.
\bid{doi={10.1287/opre.1040.0152}, issn={0030-364X}, mr={2104141}}
\end{barticle}
%
\bptok{imsref}%
\endbibitem

\bibitem{stotez}
%
\begin{barticle}[mr]
\bauthor{\bsnm{Stolyar},~\bfnm{Alexander~L.}\binits{A.~L.}} \AND
\bauthor{\bsnm{Tezcan},~\bfnm{Tolga}\binits{T.}}
(\byear{2010}).
\btitle{Control of systems with flexible multi-server pools: A shadow
routing approach}.
\bjournal{Queueing Syst.}
\bvolume{66}
\bpages{1--51}.
\bid{doi={10.1007/s11134-010-9183-0}, issn={0257-0130}, mr={2674107}}
\end{barticle}
%
\bptok{imsref}%
\endbibitem

\bibitem{wil}
%
\begin{bincollection}[mr]
\bauthor{\bsnm{Williams},~\bfnm{R.~J.}\binits{R.~J.}}
(\byear{2000}).
\btitle{On dynamic scheduling of a parallel server system with complete
resource pooling}.
In \bbooktitle{Analysis of Communication Networks: Call Centres,
Traffic and Performance ({T}oronto, {ON}, 1998)}.
\bseries{Fields Inst. Commun.}
\bvolume{28}
\bpages{49--71}.
\bpublisher{Amer. Math. Soc.},
\blocation{Providence, RI}.
\bid{mr={1788708}}
\end{bincollection}
%
\bptok{imsref}%
\endbibitem

\end{thebibliography}
\end{document}